%Paper: funct-an/9512004
%From: arveson@math.berkeley.edu (Bill Arveson)
%Date: Tue, 19 Dec 1995 16:23:36 -0800 (PST)

%NOTE: This document requires AMSTeX 2.0, the
%accompanying version of the amsppt style macros,
%and the font package amsfonts 2.0. If you try
%to typeset it using an older version of TeX, it
%probably won't work.
%
%%%%%%%%%%%%%%%%%%%  BEGINNING OF min.tex  %%%%%%%%%%%
%
\input amstex
\documentstyle{amsppt}
\loadbold
\def\cstar{$C^*$-algebra}
\def\esg{$E_0$-semigroup}

\def\<{\left<}										%for inner products
\def\>{\right>}

\magnification=\magstep 1

\topmatter
\title Minimal $E_0$-semigroups
\endtitle

\author William Arveson
\endauthor

\affil Department of Mathematics\\
University of California\\Berkeley CA 94720, USA
\endaffil

\date 6 December 1995
\enddate
\thanks This research was supported by
NSF grants DMS92-43893 and DMS95-00291
\endthanks
\abstract
It is known that every semigroup of normal completely
positive maps of a von Neumann can be ``dilated"
in a particular way to
an \esg\ acting on a larger von Neumann algebra.
The \esg\ is not uniquely determined by the completely
positive semigroup; however, it is unique (up to conjugacy)
provided that certain conditions of {\it minimality}
are met.  Minimality is a subtle property, and it is often
not obvious if it is satisfied
for specific examples even in the simplest case where
the von Neumann algebra is $\Cal B(H)$.

In this paper we clarify these
issues by giving a new characterization of
minimality in terms projective cocycles and their limits.
Our results are valid for semigroups of endomorphisms
acting on arbitrary von Neumann algebras with separable predual.
\endabstract

\endtopmatter
\vfill\eject
%Replace \pagebreak below with the line above
%to fill the lower part of the title page with
%space, rather than stretching it.
%\pagebreak

\document

\subheading{1.  Dilations and compressions}
Let $M$ be a von Neumann algebra with separable predual.  When
it is convenient to do so, we will consider that $M$ is a
von Neumann subalgebra of the algebra $\Cal B(H)$ of all bounded
operators on a separable Hilbert space $H$ which
contains the identity operator $\bold 1$.  The separability of
$H$ will be essential for some of the results below.

By an \esg\ acting on $M$ we mean a family of normal
$*$-endomorphisms $\alpha = \{\alpha_t: t\geq 0\}$ of $M$
satisfying $\alpha_t(\bold 1) = \bold 1$ for every $t\geq 0$,
which obeys the semigroup property $\alpha_{s+t} = \alpha_s\alpha_t$,
and which is continuous in the sense
that for every $a\in M$ and every pair of vectors $\xi,\eta\in H$,
the function $t\in [0,\infty)\mapsto \<\alpha_t(a)\xi,\eta\>$
is continuous.  We will also consider
semigroups $\phi = \{\phi_t: t\geq0\}$ of
normal completely positive maps acting on certain
von Neumann subalgebras $N\subseteq M$.  These subalgebras
will normally not contain the unit of $M$; but we will
require that such semigroups $\phi$ be unital in the sense that
$\phi_t(\bold 1_N) = \bold 1_N$, $t\geq 0$,
and that they should satisfy the natural
continuity property cited above.
We will refer to such a semigroup $\phi_t: N\to N$
simply as a {\it completely positive semigroup}.

$\alpha$ can be compressed to {\it certain} hereditary subalgebras
of $M$ so as to give a completely positive semigroup
as follows.  Let $M_0 = pMp$ be a hereditary
von Neumann subalgebra of $M$ with unit $p$.  The natural projection
$$
E_0: M\to M_0
$$
of $M$ onto $M_0$ is defined by $E_0(a) = pap$.  $E_0$ carries the
unit of $M$ to that of $M_0$, and we have
$$
E_0(axb) = aE_0(x)b,\qquad a,b\in M_0, x\in M.
$$
Fix $t\geq 0$.  It is an elementary exercise to show that in order
for there to exist a linear map $\phi_t: M_0\to M_0$ satisfying
$E_0\circ \alpha_t = \phi_t\circ E_0$ it is necessary and sufficient
that $\alpha_t(p)\geq p$.  Thus we will be concerned with
hereditary subalgebras $M_0 = pMp$ for which the projection
$p$ is {\it increasing} in the sense that
$$
\alpha_t(p) \geq p,\qquad t\geq 0.  \tag{1.1}
$$

In this case one can define a family $\phi = \{\phi_t: t\geq 0\}$ of
completely positive maps of $M_0$ by compressing each map
$\alpha_t$ to $M_0$,
$$
\phi_t(a) = p\alpha_t(a)p,\qquad a\in M_0,t\geq 0.  \tag{1.2}
$$
Since $\phi_t(p) = p$, we may consider
$\phi_t$ to be a unital map of $M_0$.  Note too
that the family $\phi$ has the semigroup
property $\phi_{s+t}=\phi_s\phi_t$ for $s,t\geq 0$.  Indeed,
since $p\alpha_s(p) = \alpha_s(p)p = p$ by
(1.1) we have
$$
\phi_s\phi_t(a) =
p\alpha_s(p\alpha_t(a)p)p
=p\alpha_s(p)\alpha_{s+t}(a)\alpha_s(p)p = \phi_{s+t}(a),
$$
for every $a\in M_0$.  Thus $\phi$ is a completely positive
semigroup acting on $M_0$.  Throughout this paper we will be
concerned with properties of completely positive semigroups which
can be obtained from a fixed \esg\ in this particular way.

\proclaim{Definition 1.3}
Let $p$ be a projection in $M$ satisfying (1.1), and let
$pMp$ be the corresponding hereditary subalgebra.
The completely positive semigroup $\phi = \{\phi_t: t\geq 0\}$
defined on $pMp$ by (1.2) is called a {\bf compression} of $\alpha$,
and $\alpha$ is called a {\bf dilation} of $\phi$.
\endproclaim

We emphasize that the notion of a compression of $\alpha$
to a subalgebra has meaning only when the subalgebra is
(1) a hereditary subalgebra $pMp$ of $M$ and (2) $p$ is a
projection satisfying (1.1).  It is possible to conjure up
other definitions of ``compressions" of $\alpha$
and ``dilations" of $\phi$.  For example,
one might imagine using a conditional expectation from $M$
onto a {\it unital} subalgebra $N$ of $M$ to attempt to
define a semigroup of completely positive maps of $N$
under appropriate conditions.
While there has been some limited success for such endeavors
\cite{8},\cite{9}, a recent theorem of B. V. R. Bhat \cite{5}
has led us to the conclusion that the proper context for this
kind of dilation theory is the context of Definition 1.3.

More precisely, Bhat has shown that for every
completely positive semigroup $\phi = \{\phi_t: t\geq 0\}$
acting on a von Neumann algebra $M_0$, there is a larger von Neumann
algebra $M$ containing $M_0$ as a {\it hereditary} subalgebra
$M_0 = pMp$ and an \esg\ $\alpha = \{\alpha_t: t\geq 0\}$ acting
on $M$ which satisfies $\alpha_t(p)\geq p$ for every $t$
and is such that
$\phi$ is obtained from $\alpha$ by compression.
Thus we have a dilation theory which resembles
the more familiar dilation theory for operator semigroups which
asserts that every semigroup of contraction operators
on a Hilbert space can be dilated to a semigroup of isometries on a
larger Hilbert space.

In the case of operator semigroups there is
a simple notion of {\it minimal} isometric dilation, and two minimal
isometric dilations of the same contraction semigroup are naturally
unitarily equivalent.
There is an analogous notion of minimality
in the current setting and there is an analogous
uniqueness result for the minimal \esg\ dilation of a
completely positive semigroup \cite{5}, \cite{6}.
However, these considerations for completely positive semigroups
and their \esg\ dilations are much more subtle than their
counterparts in operator theory.  To illustrate the level
of subtlety, recall that a semigroup of contraction operators
can actually be dilated further to a semigroup of unitary operators.
That is because any semigroup of isometries acting on a Hilbert
space is the restriction to an invariant subspace of a
semigroup of unitary operators acting on a larger Hilbert
space.  The unitary semigroup is unique (up to a natural
unitary equivalence) if the invariant subspace is ``minimal".
Nothing like that is true in this setting.  Indeed,
an \esg\ acting on a type $I$ factor $M$ does not have
a {\it natural} extension to a semigroup of
automorphisms of a larger type $I$ factor which contains
$M$ as a unital subfactor.  It is true
that every \esg\ can be so extended, but the construction
of the extension is quite indirect and there is
apparently no uniqueness of such extensions (see
\cite{2}, or see \cite{4} for a different proof).
Finally, the notion
of minimality introduced in \cite{5}, \cite{6} is geared to
quantum probability theory and does not lend itself
readily to the natural questions that arise in the theory of \esg s.

The purpose of this paper is to clarify the issue of minimality
for \esg s acting on arbitrary von Neumann algebras, and to
give a new characterization of minimality in terms
of the natural objects of operator algebras.

Let $p$ be an increasing projection in $M$ and let
$\phi = \{\phi_t: t\geq 0\}$ be the compression of $\alpha$
to $pMp$.  Notice that $\phi$ is itself an \esg\ (acting on
$pMp$) iff for every $t\geq 0$ we have
$$
\phi_t(ab) = \phi_t(a)\phi_t(b), \qquad a,b \in pMp.  \tag{1.4}
$$

\proclaim{Definition 1.5}
Let $\alpha$ be an \esg\ acting on a von Neumann algebra $M$.
A compression $\phi$ of $\alpha$ to a hereditary subalgebra
which satisfies property (1.4) is called {\bf multiplicative}.
\endproclaim

Suppose that $\phi$ is a compression of an \esg\ $\alpha$ to
a hereditary subalgebra $pMp$ of $M$, and that $q$ is an increasing
projection such that $q\geq p$ and the compression of $\alpha$
to $qMq$ is {\it multiplicative}.  Then may consider that the
compression of $\alpha$ to the intermediate subalgebra $qMq$
is itself an \esg\ which has $\phi$ as a compression.

\proclaim{Definition 1.6}
Let $\alpha$ be an \esg\ acting on a von Neumann algebra
$M$, let $p$ be an increasing projection in $M$, and
let $\phi$ be the completely positive semigroup on
$pMp$ obtained by compression.  $\alpha$ is said to be
{\bf minimal} over $\phi$ if the only increasing
projection $q\in M$ which satisfies $q\geq p$, and is such that
the comression of $\alpha$ to $qMq$ is multiplicative,
is the projection $q=\bold1$.
\endproclaim

In order to discuss minimality further, one needs to know
more about increasing projections which define multiplicative
compressions.  Here is the simplest class of examples.
Let $p$ be a projection
of $M$ which is {\it fixed} under $\alpha$ in the sense that
$$
\alpha_t(p)=p,\qquad t\geq 0.
$$
In this case it is clear that the compression of $\alpha$ to the
hereditary subalgebra $pMp$ is multiplicative.
However, such projections $p$ do not exhaust the
possibilities as the following observation shows.

\proclaim{Proposition 1.7}
Let $p$ be an increasing projection in $M$ and let $\phi$
be the compression of $\alpha$ to $pMp$.  Then $\phi$ is
multiplicative iff $p$ commutes with $\alpha_t(pMp)$
for every $t\geq 0$.
\endproclaim
\demo{proof}
Let $\phi_t(a) = p\alpha_t(a)p$, $a\in pMp$.  If $p$
commutes with $\alpha_t(pMp)$ then it is clear that (1.4)
is satisfied.  Conversely, if (1.4) is satisfied then for
every $a\in pMp$ we have
$$
p\alpha_t(a)^*(\bold 1-p)\alpha_t(a)p = p\alpha_t(a^*a)p
-p\alpha_t(a^*)p\alpha_t(a)p = \phi_t(a^*a)-\phi_t(a^*)\phi_t(a)=0,
$$
and hence $(\bold 1-p)\alpha_t(a)p = 0$.  Thus the range of $p$
is invariant under the self-adjoint family of operators
$\alpha_t(pMp)$, hence $p\in \alpha_t(pMp)^\prime$ \qed
\enddemo

While the criterion of Proposition 1.7 is quite specific,
it does not provide useful information for finding
the multiplicative compressions of $\alpha$.
Notice for example that the family
of von Neumann algebras $\alpha_t(pMp)$ appearing
there neither increases nor decreases
with $t$, because while the projections $\alpha_t(p)$
increase with $t$ the von Neumann algebras
$\alpha_t(M)$ {\it decrease}
with $t$.  In particular, Proposition 1.7 provides
no insight into the order structure of the family
of multiplicative compressions of $\alpha$.
In section 3 we will prove the following
two results concerning minimality.

\proclaim{Theorem A}
Let $\alpha$ be an \esg\ acting on a von Neumann algebra
$M$ with separable predual, let $p$ be an increasing
projection in $M$, and let $\phi$ be the compression of $\alpha$
to the hereditary subalgebra $pMp$.

There is an increasing projection $p_+\geq p$
which defines a multiplicative compression
of $\alpha$ such that if $q$ is any other increasing
projection satisfying $q\geq p$ for which the compression
of $\alpha$ to $qMq$ is multiplicative, then $q\geq p_+$.

The compression of $\alpha$ to $p_+Mp_+$
defines an \esg\ dilation of $\phi$ which is minimal over $\phi$.
\endproclaim

\remark{Remark}Consider $M$ to be a subalgebra
of $\Cal B(H)$.
The following result identifies the subspace $p_+H$
in concrete terms and gives an algebraic criterion
for minimality in the important case where $M$ is
a factor.
\endremark

\proclaim{Theorem B}
Let $\alpha$ be an \esg\ acting on a factor $M\subseteq \Cal B(H)$,
$H$ being a separable Hilbert space.  Let $p\in M$ be an increasing
projection and let $\phi$ be the compression of $\alpha$ to $pMp$.
The following are equivalent.
\roster
\item
$\alpha$ is minimal over $\phi$.
\item
$H$ is spanned by
$$
\{\alpha_{t_1}(a_1)\alpha_{t_2}(a_2)\dots
\alpha_{t_n}(a_n)\xi:  a_1,\dots, a_n\in pMp,
t_k\geq 0, n\geq 1, \xi\in pH \}.
$$
\item
$M$ is generated as a von Neumann algebra by the set of operators
$$
\{\alpha_t(a): a\in pMp, t\geq 0\}.
$$
\endroster
\endproclaim

\subheading{2.  Projective cocycles}
Theorems A and B will be proved in section 3.
They depend on properties of certain families
of projections satisfying a cocycle equation.

\proclaim{Definition 2.1}Let $\alpha$ be an \esg\ acting on $M$.
A {\bf projective cocycle} is a family of nonzero
projections $\{p_t: t>0\}$
in $M$ satisfying the following two conditions
$$
\align
p_t &\in \alpha_t(M)^\prime \tag{2.1.1}\\
p_{s+t} &= p_s\alpha_s(p_t),\qquad s,t> 0.  \tag{2.1.2}
\endalign
$$
\endproclaim

Notice that we have imposed no regularity condition on
the behavior of $p_t$ with respect to $t$.  This will give
us the flexibility we need for constructing examples.  Nevertheless,
projective cocycles are continuous:

\proclaim{Proposition 2.2}
Let $p = \{p_t: t>0\}$ be a projective cocycle.  Then
$p_t$ is a strongly continuous function of $t\in (0,\infty)$,
and $p_t$ tends strongly to $\bold 1$ as $t\to 0+$.
\endproclaim

\demo{proof}
The family of projections $p = \{p_t: t>0\}$
determines a family of normal self-adjoint maps
$\beta = \{\beta_t: t>0\}$ of $M$ into itself by way of
$$
\beta_t(a) = p_t\alpha_t(a),\qquad a\in M, t>0.
$$
Because of (2.1.1) each $\beta_t$ is an endomorphism
of $M$, and (2.1.2) implies that $\beta$ has the semigroup
property $\beta_{s+t} = \beta_s\beta_t$, $s,t>0$.  We have
$\beta_t(\bold 1) = p_t$ for every $t>0$.

Notice next that for fixed $\xi,\eta\in H$, the function
$t\in (0,\infty)\mapsto \<p_t\xi,\eta\>$ is
Borel-measurable.  Indeed, because of (2.1.1) and (2.1.2),
$$
p_{s+t} = p_s\alpha_s(p_t) \leq p_s
$$
for every $s,t >0$ and hence $p_t$ is decreasing in $t$.
It follows that for every $\xi \in H$, the function
$$
t\in (0,\infty) \mapsto \<p_t\xi,\xi\>\in \Bbb R^+
$$
is decreasing, therefore continuous except on a countable set,
therefore measurable.  The assertion about measurability of
$t\mapsto \<p_t\xi,\eta\>$ follows by polarization.

This implies that $\beta_t$ is weakly measurable in $t$
in the sense that
for every $\xi,\eta\in H$ and every $a\in M$, the function
$\<\beta_t(a)\xi,\eta\>$ is Borel measurable.  It follows that
for every normal linear functional $\rho\in M_*$,
$$
t\in (0,\infty) \mapsto \rho(\beta_t(a))\in \Bbb C
$$
is measurable.  Since $M_*$ is separable
we can apply Proposition 2.5 of
\cite{1} to conclude that for every $\rho\in M_*$ we have
$$
\lim_{t\to 0+}\|\rho\circ\beta_t - \rho\| = 0,
$$
and
$$
\lim_{t\to t_0}\|\rho\circ\beta_t - \rho\circ\beta_{t_0}\| = 0
\qquad \text{ for every }t_0>0.
$$
In particular, taking $\rho(a) = \<a\xi,\eta\>$ for
fixed $\xi,\eta\in H$  we
conclude that the function
$\<p_t\xi,\eta\>=\<\beta_t(\bold 1)\xi,\eta\>$ is continuous
in $t$ on the interval $(0,\infty)$ and that it
tends to $\<\xi,\eta\>=\<\bold 1\xi,\eta\>$ as $t\to 0+$.

The strong continuity of $\{p_t\}$ asserted in Prop. (2.2)
follows because the strong and weak operator topologies coincide on the
set of projections\qed
\enddemo

Thus, one may always assume that a projective cocycle
$p = \{p_t: t\geq 0\}$ is defined on the entire nonnegative
real axis, and satisfies the following two conditions
in addition to the two properties of Definition 2.1:
$$
\align
p_0 &= \bold 1,\tag{2.1.3}\\
t\in [0,\infty) &\mapsto p_t \text{ is strongly continuous.}  \tag{2.1.4}
\endalign
$$

\remark{Remarks}
Such cocycles have arisen in Powers' recent work \cite{12}
on semigroups of endomorphisms of type $I$ factors $M$.  Given
such a cocycle $p=\{p_t: t\geq 0\}$, one can form
the associated semigroup
of (nonunital) endomorphisms of $M$
$$
\beta_t(a) = p_t\alpha_t(a).
$$
Powers calls such a semigroup a {\it compression} of $\alpha$, and he
shows that the set of all compressions of $\alpha$ is a conditionally
complete lattice with respect to its natural ordering.  The compressions
of particular interest in \cite{12} are the {\it minimal} ones,
i.e., those of the form
$$
\beta_t(a) = u_tau_t^*
$$
where $\{u_t:t\geq 0\}$ is a semigroup of isometries satisfying
$$
\alpha_t(a)u_t = u_ta,\qquad a\in M, t\geq 0.
$$
Notice that in this case the projective cocycle $p$ is related
to the semigroup $u$  by $p_t = u_tu_t^*$, $t\geq 0$.

We will not use Powers' terminology in this paper because we are
concerned with dilation theory, and  in dilation
theory the term {\it compression} carries a
somewhat broader meaning.
Moreover, our need for projective cocycles
has grown from considerations that are
quite different from those of \cite{12}, and it will be
more convenient for us to deal directly with the cocycles
rather than with their associated semigroups.
\endremark

\proclaim{Definition 2.3}
Let $p = \{p_t: t>0\}$ and $q = \{p_t: t>0\}$ be two projective
cocycles.  We write $p\leq q$ if $p_t\leq q_t$ for every $t > 0$.
\endproclaim

The following proposition gives a general procedure for constructing
projective cocycles in arbitrary von Neumann algebras
from families of projections having somewhat less structure.

\proclaim{Proposition 2.4}Let $\alpha$ be an \esg\ acting on
a von Neumann algebra $M$ and let
$\{f_t: t>0\}$ be a family of nonzero projections in $M$
satisfying
$$
\align
f_t &\in \alpha_t(M)^\prime \tag{2.4.1}\\
f_{s+t} &\leq f_s\alpha_s(f_t),\qquad s,t>0.  \tag{2.4.2}
\endalign
$$
Fix $t>0$ and consider the set of all finite partitions
$$
\Cal P = \{0<t_1<t_2<\dots<t_n=t\}
$$
 of the interval $[0,t]$ as an increasing directed set
in the usual sense.  For such a partition $\Cal P$, define
an operator $f_\Cal P$ by
$$
f_\Cal P = f_{t_1}\alpha_{t_1}(f_{t_2-t_1})\alpha_{t_2}(f_{t_3-t_2})
\dots \alpha_{t_{n-1}}(f_{t_n-t_{n-1}}).  \tag{2.4.3}
$$
$f_\Cal P$ is a projection and
$\Cal P_1 \subseteq \Cal P_2 \implies f_{\Cal P_1}\leq f_{\Cal P_2}$.
Thus we can define a projection $p_t$ by
$$
p_t = \sup_\Cal P f_\Cal P = \lim_\Cal P f_\Cal P.
$$
The family $p = \{p_t: t>0\}$ is a projective cocycle, and is
the smallest projective cocycle $p$ such that
$f_t\leq p_t$ for every $t>0$.
\endproclaim

\demo{proof}
Let $s,t>0$ and let $a$ and $b$ be operators in $M$ such that $a$ commutes
with $\alpha_s(M)$ and $b$ commutes with $\alpha_t(M)$.  Then $a$
commutes with $\alpha_s(b)$ and note that
the product $a\alpha_s(b)$ commutes with
$\alpha_{s+t}(M)$.  Indeed, for arbitrary $c\in M$ we have
$$
a\alpha_s(b)\alpha_{s+t}(c) = a\alpha_s(b\alpha_t(c))
=\alpha_s(\alpha_t(c)b)a = \alpha_{s+t}(c)\alpha_s(b)a
= \alpha_{s+t}(c)a\alpha_s(b).
$$

Now fix $t>0$.
It follows from the preceding remarks
that the operator $f_\Cal P$ of (2.4.3) belongs to
$M\cap \alpha_t(M)^\prime$; moreover, the $n$ factors of $f_\Cal P$
on the right side of (2.4.3) are mutually commuting projections.
Thus $f_\Cal P$ is a projection in $M\cap \alpha_t(M)^\prime$.

To show that $f_\Cal P$ increases with $\Cal P$ it is enough
to show that if a given partition
$\Cal P = \{0<t_1<\dots<t_n=t\}$ is refined by adjoining to it
a single point $\tau$, then $f_\Cal P$ increases.
In turn, that reduces
to the following assertion.  For $k = 1,2,\dots,n$ and
$t_{k-1} < \tau < t_k$ (where $t_0$ is taken as $0$),
$$
f_{t_k-t_{k-1}} \leq f_{\tau-t_{k-1}}\alpha_{\tau-t_{k-1}}(f_{t_k-\tau}).
$$
The latter is immediate from the hypothesis (2.4.2).

Thus the net $f_\Cal P$ increases with $\Cal P$ and we
can define a projection $p_t\in M\cap\alpha_t(M)^\prime$ as asserted.
The cocycle property (2.1.2) follows immediately from the definition
of the family $\{p_t: t>0\}$.

We obviously have $f_t\leq p_t$ for every $t>0$.  Finally, suppose
$q = \{q_t: t>0\}$ is another projective cocycle satisfying
$f_t\leq q_t$ for every $t>0$.  Fix $t>0$ and let
$\Cal P = \{0<t_1<\dots<t_n=t\}$ be
a partition of the interval $[0,t]$.  Then for every $k=1,2,\dots,n$
we have $f_{t_k - t_{k-1}}\leq q_{t_k - t_{k-1}}$, and hence
$f_\Cal P\leq q_\Cal P$.  On the other hand, the cocycle property of
$q$ implies that $q_\Cal P = q_t$.  Hence $f_\Cal P\leq q_t$ and
we deduce the desired inequality
$$
p_t = \sup_\Cal P f_\Cal P \leq q_t.
$$
\qed
\enddemo

The following result will be important for section 3.
It implies that certain projections in $M$ naturally
give rise to projective cocycles.

\proclaim{Corollary 2.5}
Let $e$ be a nonzero projection in $M$ satisfying
$\alpha_t(e)\geq e$ for every $t\geq 0$.  For each
$t>0$ let $f_t$ be the smallest projection
in $M\cap\alpha_t(M)^\prime$
which dominates $e$, i.e.,
$$
f_t = [\alpha_t(a)e\xi: a\in M, \xi \in H].
$$
Then $f_{s+t}\leq f_s\alpha_s(f_t)$ for every $s,t>0$.  The
projective cocycle $p = \{p_t: t>0\}$ of Proposition 2.4
is the smallest projective cocycle satisfying $p_t\geq e$
for every $t>0$.
\endproclaim

\demo{proof}
It is obvious that $f_t$ commutes with $\alpha_t(M)$, and
the double commutant theorem implies that $f_t\in M$.  Hence
$f_t\in M\cap\alpha_t(M)^\prime$.

To see that $f_{s+t}\leq f_s\alpha_s(f_t)$, fix $s,t>0$.
By the argument at the beginning of the proof of
proposition 2.4, $f_s\alpha_s(f_t)$ is a projection
in $M\cap\alpha_{s+t}(M)^\prime$.  We claim that
$e\leq f_s\alpha_s(f_t)$.  Indeed, $e\leq f_s$ follows from
the definition of $f_s$, and since $e\leq f_t$ implies
$\alpha_s(e)\leq \alpha_s(f_t)$ we have
$e\leq \alpha_s(e)\leq \alpha_s(f_t)$.  Hence
$e\leq f_s\alpha_s(f_t)$.
Since $f_{s+t}$ is the smallest projection in
$M\cap\alpha_{s+t}(M)^\prime$ which dominates $e$ we
have the asserted inequality $f_{s+t}\leq f_s\alpha_s(f_t)$\qed
\enddemo

%\documentstyle{amsppt}
%\loadbold
%\def\cstar{$C^*$-algebra}
%\def\esg{$E_0$-semigroup}
%\def\Log{\text{Log}\,}
%\def\<{\left<}										%for inner products
%\def\>{\right>}
%\define\sumplus{\sideset \and^\oplus \to\sum}
%\magnification=\magstep 1

\subheading{3.  Minimality}

Throughout the section, $\alpha$ will denote an \esg\ acting
on a von Neumann algebra $M\subseteq\Cal B(H)$ and
$p$ will denote a projection in $M$ satisfying
$$
\alpha_t(p)\geq p,\qquad t\geq 0.  \tag{3.1}
$$
$M_0$ will denote the hereditary subalgebra $pMp$.
We will be concerned with the (perhaps nonunital) von Neumann
algebra $M_+$ generated by $M_0$ and its translates under $\alpha$,
$$
M_+ = \overline{\text{span}}\{\alpha_{t_1}(a_1)\alpha_{t_2}(a_2)\dots
\alpha_{t_n}(a_n): a_k\in M_0, t_k\geq 0, n=1,2,\dots\},
$$
the bar denoting closure in the weak operator topology.  It is
obvious that $\alpha_t(M_+)\subseteq M_+$ for every $t\geq 0$, and
the unit of $M_+$ is the projection
$$
p_\infty = \lim_{t\to\infty} \alpha_t(p).  \tag{3.3}
$$
There is a somewhat smaller projection that is
of greater importance, namely
$$
p_+ = [M_+pH] =
[\alpha_{t_1}(a_1)\alpha_{t_2}(a_2)\dots\alpha_{t_n}(a_n)\xi:
a_k\in M_0, t_k\geq 0, \xi\in pH\}.  \tag{3.4}
$$
\remark{Remarks}
$p_+$ is the unit of the two-sided ideal
$\overline{\text{span}}M_+pM_+$ of $M_+$
generated by $p$.  Thus, $p_+$ belongs to the
center of $M_+$, and in fact is the smallest
central projection $c$ in $M_+$ satisfying
$p\leq c$.

We will eventually show that $p_+$ is an
increasing projection which defines a
multiplicative compression
of $\alpha$.  Neither of these assertions
is apparent from (3.4).  We deduce these
properties from the following result which
gives a ``formula" for $p_+$.
\endremark

\proclaim{Theorem 3.5}Let $p$ be a projection in $M$ which
satisfies (3.1), let $p_\infty$ and $p_+$ be defined as in
(3.3) and (3.4) respectively.  Let $q = \{q_t: t>0\}$ be
the smallest projective cocycle satisfying $q_t\geq p$ for every
$t\geq 0$ as in Corollary 1.5, and set
$$
q_\infty = \lim_{t\to\infty}q_t.
$$
Then $p_\infty$ belongs to the tail von Neumann algebra
$M_\infty = \cap_{t\geq 0}\alpha_t(M)$, $q_\infty$
belongs to its relative commutant in $M$, and we have a factorization
$$
p_+ = p_\infty q_\infty.
$$
\endproclaim

\remark{Remarks}
Recall that since $q$ is a projective cocycle, $q_t$ must
be a decreasing function of $t$ and hence the strong limit
$q_\infty = \lim_{t\to\infty}q_t$ exists.
\endremark

\demo{proof}
Notice first that since $\alpha_t(p)\in \alpha_t(M)$ and since the
von Neumann algebras $\alpha_t(M)$ decrease as $t$ increases, it
follows that $p_\infty= \lim_{t\to\infty}p_t\in M_\infty$.
Since $q_t$ belongs to
$M\cap\alpha_t(M)^\prime \subseteq M\cap M_\infty^\prime$ for all $t$
we see that $q_\infty = \lim_{t\to\infty}q_t\in M\cap M_\infty^\prime$.
In particular, the projections $p_\infty$ and $q_\infty$ must commute.

We show first that $p_+\leq p_\infty q_\infty$.  Since $p_+\leq p_\infty$
is obvious, it suffices to show that $p_+H\leq q_\infty H$. Considering
the definition  of $p_+$ and the fact that $pH\leq q_\infty H$, it suffices
to show that the subspace $q_\infty H$ is invariant under any operator
in any one of the von Neumann algebras $\alpha_t(M_0)$, $t> 0$, i.e.,
that $q_\infty$ commutes with $\cup_{t>0}\alpha_t(M_0)$.  For
that, fix $t>0$.  If we pass $s$ to $\infty$ in the cocycle formula
$$
q_t\alpha_t(q_s) = q_{t+s}
$$
and use normality of $\alpha_t$ we obtain
$$
q_t\alpha_t(q_\infty) = q_\infty.  \tag{3.6}
$$
It follows that for any $a\in M_0$ we have
$$
\alpha_t(a)q_\infty = \alpha_t(a)q_t\alpha_t(q_\infty)
= q_t\alpha_t(a)\alpha_t(q_\infty)= q_t\alpha_t(aq_\infty), \tag{3.7}
$$
where we have used $q_t\in \alpha_t(M)^\prime$.  Now since $p\leq q_\infty$
and since $a = pap$ we have $aq_\infty = a = q_\infty a$.  Thus
we can replace the right side of (3.7) with
$$
q_t\alpha_t(q_\infty a) = q_t\alpha_t(q_\infty) \alpha_t(a)
= q_\infty \alpha_t(a).
$$
Thus $\alpha_t(a)$ commutes with $q_\infty$ as required.

It remains to show that $p_\infty q_\infty \leq p_+$.  For that,
it suffices to show that for every $t>0$ we have
$$
q_t\alpha_t(p_+)\leq p_+.  \tag{3.8}
$$
Indeed, assuming that (3.8) has been established we deduce
 $q_t\alpha_t(p) \leq p_+$ for every $t$ (because $p\leq p_+$);
noting that
$q_t \downarrow q_\infty$ and $\alpha_t(p)\uparrow p_\infty$
as $t$ increases to $+\infty$, we
may take the strong limit on $t$ in
the previous formula to obtain
the desired inequality $q_\infty p_\infty \leq p_+$.

In order to prove (3.8), we require
\proclaim{Lemma 3.9}For each $t>0$ let $f_t$ be the
projection onto the subspace
$$
[\alpha_t(a)p\xi: a\in M, \xi \in H].
$$
Then $f_t\alpha_t(p_+)\leq p_+$.
\endproclaim
\demo{proof}
We have already pointed out in the remarks following (3.4)
that $p_+$ is the unit of the ideal $\overline{\text{span}}M_+pM_+$
in $M_+$.  Thus it suffices to show that
$$
f_t\alpha_t(M_+pM_+)H\subseteq p_+H.
$$
Since $f_t$ commutes with $\alpha_t(M)$ the left side
is contained in
$$
\alpha_t(M_+pM_+)f_tH \subseteq [\alpha_t(M_+pM_+)\alpha_t(M)pH]
\subseteq [\alpha_t(M_+pM_+M)pH].
$$
Noting that $p = \alpha_t(p)p$ the latter is
$$
\align
[\alpha_t(M_+pM_+Mp)pH]&\subseteq [\alpha_t(M_+pMp)pH]
\subseteq [\alpha_t(M_+)pH]\\
&\subseteq [M_+pH] = p_+H,
\endalign
$$
as asserted\qed
\enddemo

For every $t>0$ we define a normal linear mapping
$\beta_t: M\to M$ as follows,
$$
\beta_t(a) = f_t\alpha_t(a).
$$
Since $f_t$ commutes with $\alpha_t(M)$, $\beta_t$ is an
endomorphism of the $*$-algebra structure of $M$ for
which $\beta_t(\bold 1) = f_t$, but it
is not a semigroup because $\{f_t: t>0\}$ does not satisfy
the cocycle condition (1.1.2).  However, because of Lemma 3.9
we have
$$
\beta_t(p_+) \leq p_+, \qquad t>0.
$$

Now fix $t>0$ and let $\Cal P = \{0=t_0<t_1<\dots<t_n=t\}$ be
a partition of the interval $[0,t]$.  By iterative the preceding
formula we find that
$$
\beta_{t_1}\beta_{t_2-t_1}\beta_{t_3-t_2}
\dots \beta_{t_n-t_{n-1}}(p_+) \leq p_+.
\tag{3.10}
$$
In the notation of Proposition 1.4, the left side of (3.10) is
$$
f_{t_1}\alpha_{t_1}(f_{t_2-t_1})\alpha_{t_2}(f_{t_3-t_2})
\dots \alpha_{t_{n-1}}(f_{t_n-t_{n-1}})\alpha_t(p_+)
= f_\Cal P \alpha_t(p_+).
$$
Using 1.4 and 1.5, we make take the limit on $\Cal P$ in (3.10)
obtain the required inequality (3.8)
$$
q_t\alpha_t(p_+) = \lim_\Cal P f_\Cal P\alpha_t(p_+)\leq p_+,
$$
completing the proof of Theorem 3.5\qed
\enddemo

We can now deduce the following result, which paraphrases
Theorem A from section 1.

\proclaim{Theorem A}Let $p$ be an increasing projection
for $\alpha$ and let $p_+\geq p$ be the projection
defined by (3.4).
$p_+$ is an increasing projection with the property
that the compression of $\alpha$ to $p_+Mp_+$ is
multiplicative.

If $r$ is another increasing projection in $M$ such that
$r\geq p$ and the compression of $\alpha$ to
$rMr$ is multiplicative, then $r\geq p_+$.
\endproclaim
\demo{proof}
Since $p_\infty = \lim_{t\to\infty}\alpha_t(p)$ is clearly
fixed under the action of $\alpha_t$ and since (3.6) implies
that $\alpha_t(q_\infty) \geq q_\infty$, we
find that $\alpha_t(p_+) = \alpha_t(p_\infty)\alpha_t(q_\infty)
\geq p_\infty q_\infty = p_+$.

To show that the compression of $\alpha$ to $p_+Mp_+$ is
multiplicative, it suffices to show that $p_+$ commutes
with $\alpha_t(p_+Mp_+)$ for every $t>0$ (Propostion 1.7).
For that, it is enough to show
that for every $a=a^*\in p_+Mp_+$ we have
$$
p_+\alpha_t(a) = q_t\alpha_t(a).  \tag{3.11}
$$
Indeed, by taking
adjoints in (3.11) we find that $\alpha_t(a)q_t = \alpha_t(a)p_+$,
and since $q$ is a projective cocycle $q_t$ must commute
with $\alpha_t(M)$.  Thus
$$
p_+\alpha_t(a) = q_t\alpha_t(a)=\alpha_t(a)q_t = \alpha_t(a)p_+,
$$
and thus $p_+\in \alpha_t(p_+Mp_+)^\prime$.

To prove (3.11), we write $p_+ = p_\infty q_\infty = q_\infty p_\infty$
and use $\alpha_t(p_\infty) = p_\infty$ to obtain
$$
p_+\alpha_t(a) = q_\infty p_\infty \alpha(t(a) =
q_\infty\alpha_t(p_\infty a) = q_\infty\alpha_t(a),\tag{3.12}
$$
because $p_\infty a = a$ for every operator $a$ in
$p_+ Mp_+\subseteq p_\infty M p_\infty$.  Using (3.6) on the
last term of (3.12) we have
$$
 q_\infty\alpha_t(a) = q_t\alpha_t(p_\infty)\alpha_t(a)
= q_t\alpha_t(q_\infty a) = q_t\alpha_t(a),
$$
since $q_\infty a = a$
for every $a\in p_+ Mp_+\subseteq q_\infty M q_\infty$.
Formula (3.11) follows.

Finally, let $r\geq p$ be another increasing projection with
the property that the compression of $\alpha$ to $rMr$ is
multiplicative.  We have to show that $p_+H \subseteq rH$.
Because of formula (3.4) for $p_+H$,
together with the fact that
$pH\subseteq rH$, it is enough to show that
$rH$ is invariant under any operator of the form
$\alpha_t(a)$ with $a\in pMp$ and $t>0$.  But
for each $t>0$,
Proposition 1.7 implies that $r$ commutes
with the set of operators $\alpha_t(rMr)$, and therefore
since $p\leq r$ we have
$$
\alpha_t(pMp)rH\subseteq \alpha_t(rMr)rH\subseteq rH,
$$
as required \qed
\enddemo

As another consequence of Theorem 3.5 we have the following
characterization of minimality in terms of projective cocycles.

\proclaim{Corollary 3.13}
Let $p$ be an increasing projection and let $\phi$ be the
compression of $\alpha$ to $pMp$.  Then $\alpha$ is minimal
over $\phi$ iff $\lim_{t\to\infty}\alpha_t(p) = \bold 1$ and
the only projective cocycle $q = \{q_t: t>0\}$ satisfying
$q_t\geq p$ for every $t>0$ is the trivial cocycle $q_t = \bold 1$.
\endproclaim
\demo{proof}
The minimality assertion is that $p_+ = \bold 1$ and from
Theorem 3.5 we have $p_+=p_\infty q_\infty$.  Thus $\alpha$
is minimal iff $p_\infty = q_\infty = \bold 1$\qed
\enddemo

\proclaim{Proposition 3.14}  Let $M_+$ be the following
von Neumann subalgebra of $M$
$$
M_+ = \overline{span}\{\alpha_{t_1}(a_1)\alpha_{t_2}(a_2)
\dots\alpha_{t_n}(a_n): a_1,\dots,a_n\in pMp, t_1,\dots,t_n\geq 0, n\geq 1\},
$$
and let $p_+$ be the projection of (3.4).  Then $p_+$
is the smallest projection in the center of $M_+$
which dominates $p$, and we have $p_+Mp_+ = M_+p_+$.
\endproclaim

\demo{proof}  While the first assertion is very
elementary, we include a proof for completeness.
Clearly $p_+\in M_+$, and since $p_+H = [M_+pH]$ is invariant
under $M_+$ we have $p_+\in M_+^\prime$.  Hence $p_+$ belongs
to the center of $M_+$ and $p\leq p_+$.  If $c$ is another
central projection in $M_+$ for which $c\geq p$, then $cH$
clearly contains $[M_+pH] = p_+H$ and hence $c\geq p_+$.

Let $R$ be the weakly closed subspace of $M$ generated by the
set of operators
$$
\{apb: a\in M_+, b\in M\}.
$$
R is a right ideal in $M$, and the range projection of $R$ is
$$
[RH] = [apH: a\in M_+] = p_+H.
$$
Hence $R = p_+M$.  Thus $p_+Mp_+$ is spanned by $RR^*$, i.e.,
$$
p_+Mp_+ = \overline{span}(M_+\cdot pMp \cdot M_+)
= \overline{span}(M_+ p M_+).
$$
The right side of the preceding formula is the two-sided ideal
in $M_+$ generated by $p$ which, by the preceding paragraph,
is $M_+p_+$\qed
\enddemo

%We now turn to the proof of
\proclaim{Theorem B}
Let $\alpha$ be an \esg\ acting on a factor $M\subseteq \Cal B(H)$,
$H$ being a separable Hilbert space.  Let $p\in M$ be an increasing
projection and let $\phi$ be the compression of $\alpha$ to $pMp$.
The following are equivalent.
\roster
\item
$\alpha$ is minimal over $\phi$.
\item
$H$ is spanned by
$$
\{\alpha_{t_1}(a_1)\alpha_{t_2}(a_2)\dots
\alpha_{t_n}(a_n)\xi:  a_1,\dots, a_n\in pMp,
t_k\geq 0, n\geq 1, \xi\in pH\}.
$$
\item
$M$ is generated as a von Neumann algebra by the set of operators
$$
\{\alpha_t(a): a\in pMp, t\geq 0\}.
$$
\endroster
\endproclaim

\demo{proof of (1)$\implies$(3)}
If $\alpha$ is minimal over $\phi$ then $p_+ = \bold 1$, and
thus by Proposition 3.14 we find that
$M = p_+Mp_+ = M_+p_+ = M_+$, hence (3).
\enddemo

\demo{proof of (3)$\implies$(2)}
Since $M = M_+$ we have $[M_+pH] = [MpH]$.  The projection
on the subspace on the right is the central carrier of
$p$, which must be $\bold 1$ because $M$ is a factor.
Therefore $p_+H = H$, as asserted in (2).
\enddemo

\demo{proof of (2)$\implies$(1)}
The hypothesis is that $p_+ = \bold 1$ which, by Theorem A,
implies that $\alpha$ is minimal over $\phi$ \qed
\enddemo

\vfill
\pagebreak

\end
\Refs
%
\ref\no 1\by Arveson, W.\paper Continuous analogues of Fock space
\jour Memoirs Amer. Math. Soc.\vol 80 no. 3\yr 1989
\endref

\ref\no 2\bysame\paper Continuous analogues of Fock space IV:
essential states\jour Acta Math.\vol 164\yr 1990\pages 265--300
\endref

\ref\no 3\bysame\paper Non-commutative flows I: dynamical
principles\paperinfo preprint November 1995
\endref

\ref\no 4\by Arveson, W. and Kishimoto, A.\paper A note on extensions
of semigroups of $*$-endomorphisms\jour Proc. A. M. S.\vol 116, no 3
\yr 1992\pages 769--774
\endref

\ref\no 5\by Bhat, B. V. R. \paper Minimal dilations of
quantum dynamical semigroups to semigroups of endomorphisms of
\cstar s \jour Trans. A.M.S. \toappear
\endref

\ref\no 6\bysame \paper On minimality of Evans-Hudson flows
\jour (preprint)
\endref

\ref\no 7\by Hudson, R. L. and Parthasarathy, K. R. \paper
Stochastic dilations of uniformly continuous completely positive
semigroups \jour Acta Appl. Math. \vol 2\pages 353--378\yr 1984
\endref

\ref\no 8\by K\"ummerer, B. \paper
Markov dilations on $W^*$-algebras \jour J. Funct. Anal.
 \vol 63\pages 139--177\yr 1985
\endref

\ref\no 9\bysame\paper Survey on a theory of non-commutative
stationary Markov processes\inbook Quantum Probability and
Applications III\publ Springer Lecture notes in Mathematics
\vol 1303\yr 1987\pages 154--182
\endref

\ref\no 10\by Mohari, A., Sinha, Kalyan B. \paper Stochastic
dilation of minimal quantum dynamical semigroups \jour
Proc. Ind. Acad. Sci. \vol 102\yr 1992\pages 159--173
\endref

\ref\no 11\by Parthasarathy, K. R. \book An introduction to quantum
stochastic calculus \publ Birkh\"auser Verlag, Basel\yr 1991
\endref

\ref\no 12\by Powers, R. T. \paper New examples of continuous spatial
semigroups of endomorphisms of $\Cal B(H)$
\paperinfo preprint 1994
\endref

\endRefs
